 \documentclass[10pt,twocolumn,twoside]{IEEEtran}
\usepackage{color}
\usepackage{graphicx,subfigure}
\usepackage{epsfig}
\usepackage{parskip}
\usepackage{enumitem}
\usepackage{times}
\usepackage{url}
\usepackage[cmex10]{amsmath}
\usepackage{amsthm,amssymb}
\usepackage{algorithm}
\usepackage{algpseudocode}

\newcommand{\comment}[1]{}

\newcommand{\Lf}{\ensuremath{L_{ff}}}
\newcommand{\Ll}{\ensuremath{L_{ll}}}
\newcommand{\Lfl}{\ensuremath{L_{fl}}}
\newcommand{\Llf}{\ensuremath{L_{lf}}}

\newcommand{\xvec}{\ensuremath{\textbf{x}}}
\newcommand{\xl}{\ensuremath{\textbf{x}_l}}
\newcommand{\xf}{\ensuremath{\textbf{x}_f}}
\newcommand{\dvec}{\ensuremath{\textbf{d}}}

\newcommand{\trace}[1]{\ensuremath{\textbf{tr}\left(#1\right)}}
\newcommand{\tp}{\ensuremath{^{\mathsf{T}}}}

\newcommand{\expec}[1]{\ensuremath{\textbf{E}\left(#1\right)}}

\newcommand{\cC}{\ensuremath{{\mathcal C}}}

\newcommand{\M}{\ensuremath{M}}

\newcommand{\mR}{\ensuremath{\mathbf{R}}}
\newcommand{\cstar}{\ensuremath{c^{\star}}}

\newcommand{\lev}[1]{\ensuremath{\text{lev}(#1)}}

\newcommand{\xbar}{\ensuremath{\overline{x}}}
\newcommand{\cN}{\ensuremath{{\mathcal N}}}

\newcommand{\Rnf}{R_{NF}}
\newcommand{\Rnc}{R_{NC}}
\newcommand{\sbar}{\overline{s}}
\newcommand{\rbar}{\overline{r}}
\newcommand{\Gbar}{\overline{G}}
\newcommand{\Ebar}{\overline{E}}
\newcommand{\Vbar}{\overline{V}}
\newcommand{\Wbar}{\overline{W}}

\newtheorem{theorem}{Theorem}
\newtheorem{lemma}[theorem]{Lemma}
\newtheorem{proposition}[theorem]{Proposition}

\newtheorem*{theor-non}{Theorem}
\newtheorem*{propo-non}{Proposition}

\begin{document}
%

\title{A Resistance Distance-Based Approach for Optimal Leader Selection in Noisy Consensus Networks}

\author{Stacy~Patterson,~\IEEEmembership{Member,~IEEE}, 
       ~Yuhao~Yi, and Zhongzhi~Zhang%
\thanks{S.~Patterson$^*$ is with the Department of Computer Science,
Rensselaer Polytechnic Institute, Troy, NY 12180. Email: \emph{sep@cs.rpi.edu}, Phone: 518-276-2054 ($^*$corresponding author)}
\thanks{Y.~Yi and Z.~Zhang are with the Shanghai Key Laboratory of Intelligent Information Processing, School of Computer Science, Fudan University, Shanghai 200433, China.  
Email: \emph{15110240008@fudan.edu.cn}, \emph{zhangzz@fudan.edu.cn}.}}

\maketitle


\begin{abstract}
We study the performance of leader-follower noisy consensus networks, and in particular, the relationship between this performance and the locations of the leader nodes.
Two types of dynamics are considered
(1) noise-free leaders, in which  leaders dictate the trajectory exactly and  followers are subject to external disturbances, and (2) noise-corrupted leaders, in which both leaders and followers are subject to external perturbations.
We measure the performance of a network by its \emph{coherence}, an $H_2$ norm that quantifies how closely the followers track the leaders' trajectory.  For both dynamics, we show a relationship between the coherence and
resistance distances in an a electrical network. Using this relationship, we derive closed-form expressions for coherence as a function of the locations of the leaders.
Further, we give analytical solutions to the optimal leader selection problem for several special classes of graphs.
\end{abstract}


\section{Introduction}
Consensus problems are an important class of problems in networked and multi-agent systems.  The consensus model has been used to study a wide range of applications, including opinion dynamics in social networks~\cite{degroot1974reaching}, information fusion in sensor networks~\cite{xiao2005scheme}, formation control~\cite{BJMP12}, and load balancing in distributed computing systems~\cite{Cybenko}. Over the past decades, much research effort has been devoted to analysis of the convergence behavior and robustness of consensus networks and to the derivation of relationships between system performance and graph theoretic properties.

A type of consensus problem that has received attention in recent years is leader-follower consensus~\cite{PB10,CABP14,LFJ14,P15,L16,7277027,PMD16}. In leader-follower systems, 
 a subset of nodes are  leaders that track  an external signal.  The leaders, in essence, dictate the desired trajectory of the network.
The remaining nodes are followers that update their states based on relative information exchanges with neighbors.  
Leader-follower dynamics can be used to model formation control where, due to bandwidth limitations, only a small subset of agents can be controlled by a system operator~\cite{BH06}. In addition, leader-follower systems can also be used to model agreement dynamics in social networks in which some subset of participants exhibit degrees of stubbornness~\cite{6732931}.  Leader-follower dynamics have also been applied to the problem of distributed sensor localization~\cite{BH09}.
In leader-follower systems, the system performance depends on the network topology and the locations of the leaders.  This dependence naturally leads to the question of how to select the leaders so as to optimize performance for a given topology.

We study the the performance of leader-follower networks where nodes are governed by consensus dynamics and are also subject to stochastic external disturbances.   We consider two types of dynamics. In the first,  referred to as \emph{noise-free leaders}, leaders are not subject to disturbances and thus track the external signal exactly. In the second dynamics, called \emph{noise-corrupted leaders}, all nodes are subject to the external perturbations. As in many works on noisy consensus networks~\cite{PB10,CBP14,PMD16,7277027}, we quantify the system performance by an $H_2$ norm that captures the steady-state variance of the node states. We call this the \emph{coherence} of the network. Coherence is related to the spectrum of the Laplacian matrix of the network; however, it is not always straightforward to relate this spectrum to the network topology and locations of leaders.

In this work, we  develop relationships between the steady-state variance for a given leader set and resistance distances in a corresponding electrical network.  A similar approach was used to study the performance of a single noise-free leader~\cite{BH06}; we generalize this notion to an arbitrary number of noise-free leaders. Further, we develop a novel resistance-distance based approach to study coherence in networks with an arbitrary number of noise-corrupted leaders.  
We use this resistance distance-based approach to  analyze the coherence for different network topologies based on resistance distances.
In special classes of graphs, we can relate the resistance distance to graph distance, which gives us the optimal leader locations in terms of the graph distances between leaders.
We also e derive closed form-expressions for the optimal  single noise-free and noise-corrupted leaders in weighted graphs, the optimal $k$ noise-free leaders in cycles and paths, the optimal two noise-free leaders in trees, and the optimal twof noise-corrupted leaders in cycles.

The leader selection problem for noise-free leaders was first posed in~\cite{PB10}. This problem can be solved by an exhaustive search over all subsets of nodes of size $k$, but this proves computationally intractable for large graphs and large $k$. 
Several works have proposed polynomial-time approximation algorithms for the $k$-leader selection problem in noise-free leader-follower systems~\cite{BH09,LFJ14,CBP14,P15}.
In particular, we note that the solution presented in \cite{CBP14} yields a leader set whose performance is within a provable bound of optimal. 
With respect to analysis for the noise-free leader selection problem, the recent work by Lin~\cite{L16} gives
asymptotic scalings of the steady-state variance in directed lattice graphs for a single noise-free leader, based on the graph distance from the leader.
Our recent work~\cite{PMD16} gives polynomial-time algorithms for \emph{optimal} $k$-leader selection in weighted, undirected cycles and path graphs.
The leader-selection problem for noise-corrupted leaders was first posed by Lin et al.~\cite{LFJ14}, who also gave heuristic-based bounds and algorithms for its solution.
 In addition, other performance measures have been considered for the leader selection problem including controllability~\cite{O15,FL16} and convergence rate~\cite{CABP14,PMD16}.

 The recent works by Fitch and Leonard~\cite{7277027,FL16} study the optimal leader selection problem for noise-free and noise-corrupted leaders.  These works also relate the steady-state variance to  a graph theoretic concept, in this case, graph centrality.  The authors define centrality measures that capture the performance of a given leader set.  They then use this analysis to identify the optimal leader sets for various classes of graphs.  We note that this work identifies the optimal single leader for noise-free and noise-corrupted graphs under slightly stronger assumptions than we make in our approach. In addition,~\cite{7277027} identifies the optimal $k$-noise free leaders in cycles, under the restriction that the number of nodes in the cycle is a multiple of $k$.  We address cycles with an arbitrary number of nodes  and provide a closed-form expression for the resulting steady-state variance for any leader set  based on the graph distance between leaders.  We view our proposed approach as complementary that in~\cite{7277027}; for some classes of networks, analysis is more straightforward under the resistance distance interpretation. Thus, our work expands the classes of networks that have known analytical solutions. 
A preliminary version of our work appeared in~\cite{P17}.  This earlier work gave analysis for noise-free leader selection in cycle and path graphs only, using the related concept of commute times of random walks rather than resistance distance. Our resistance-distance based approach greatly simplifies the analysis and presentation.

The remainder of this paper is organized as follows. Section~\ref{model.sec} describes the system model and dynamics, and it formalizes the leader selection problems.
Section~\ref{resistance.sec} describes the relationship between the system performance and  resistance distance for both noise-free and noise corrupted leaders.
This section also presents analysis of resistance distance for ``building blocks'', i.e., components of graphs, that will be used to analyze specific graph topologies.
Section~\ref{graphs.sec} gives closed-form solutions for the leader selection problem for various classes of graphs. In Section~\ref{compare.sec}, we compare the asymptotic behavior of coherence
in leader-free and leader-follower consensus networks, and in Section~\ref{treegrowing.sec}, we give an algorithm and a numerical example for increasing the size of a binary tree while maintaining the optimality of the  two noise free leaders.
Finally, we conclude in Section~\ref{conclusion.sec}.

\section{System Model and Problem Formulation} \label{model.sec}

We consider a network of $n$  agents, modeled by an undirected, connected graph ${G=(V,E,W)}$, where $V$ is the set of agents, also called nodes, and $E$ is the set of edges.
The weight of edge $(i,j)$, denoted by $w_{ij}$, corresponds to the $(i,j)^{th}$ component of the symmetric weighted adjacency matrix $W$.
We let $D$ denote the diagonal matrix of weighted node degrees, with diagonal entries $d_{ii} = \sum_{j \in V} w_{ij}$.  
The matrix $L = D-W$ is thus the weighted Laplacian matrix of the graph $G$.

 Each node $i \in V$ has a scalar-valued state $x_i$.  The objective is for all node states to track an external signal $\xbar \in \mR$.
Some subset of nodes $F \subset V$ are  \emph{followers} that update their states using noisy consensus dynamics, i.e., 
\begin{equation}
\dot{x}_{i}  = - \sum_{j \in {\cN}(i)} w_{ij} \left( x_i - x_j \right) + d_i, \label{fdyn.eq}
\end{equation}
where  $\cN_i$ denotes the neighbor set of node $i$, and $d_i$ is a zero-mean, unit variance, white stochastic noise process.
The remaining set of nodes $S = V \setminus F$ are \emph{leaders}; leader nodes have access  $\xbar$.

We write the state of the system as $\xvec \tp = [ \xl \tp~\xf \tp]$,
where $\xl$ are the leader states and $\xf$ are the follower states.  
We can then decompose the Laplacian of $G$ as:
\[
L = \left(
\begin{array}{cc}
\Ll & \Llf \\
\Lfl & \Lf
\end{array}
\right).
\]

\subsection{Noise-Free Leader Dynamics}
We consider two types of leader dynamics. In the first, called \emph{noise-free leaders}, leader states are dictated solely by $\xbar$.
  Without loss of generality,
we assume $\xbar = 0$~\cite{PB10}, so
leader nodes update their states as:
\[
\dot{x}_i = -\kappa_i (x_i - \xbar) = -\kappa_i x_i ,
\]
where $\kappa_i \in \mR^{+}$ is the weight  node $i$ gives to the external signal, sometimes referred to as the degree of stubbornness of node $i$.
The dynamics of the follower nodes can then be written as:
\begin{equation} \label{nfdyn.eq}
\dot{\textbf{x}}_f= -  \Lf \xvec  + \dvec_f,
\end{equation}
where $\Lf$ is the principle submatrix of the Laplacian corresponding to the follower nodes, and  $\dvec_f$ is the vector of noise processes for the followers.

We quantify the performance of the system for a given leader set $S$ by its coherence,
 i.e., the total steady-state variance of the follower nodes. This value  is related to $\Lf$
as follows~\cite{PB10},
\begin{equation} \label{coherence.eq}
\Rnf(S) = \lim_{t \rightarrow \infty} \sum_{i \in  (V \setminus S)} \expec{x_i(t)^2} =  \frac{1}{2} \trace{(\Lf)^{-1}}. 
\end{equation}
Note that $\Lf$ is positive definite for any $S \neq 0$~\cite{PB10}, and thus, $\Rnf(S)$ is well defined.
The total variance depends on the choice of leader nodes.  

The \emph{nose-free leader selection problem} is to identify the leader set $S$ of size at most $k$, such that 
$\Rnf(S)$ is as small as possible, i.e.,
\begin{equation}
\begin{array}{ll}
\mathrm{minimize}& \Rnf(S) \\
\mathrm{subject~to}& |S|  \leq k.
\end{array} \label{kleader.eq}
\end{equation}

\subsection{Noise-Corrupted Leader Dynamics}
We also consider dynamics with noise-corrupted leaders.  In this case, the leader nodes update their states using both consensus dynamics and the external signal, and the leader states are also subject to external disturbances. We again assume, without loss of generality, that $\xbar$ is 0. The dynamics for leader node $i$ are then:
\[
\dot{x}_i = - \sum_{j \in {\cN}(i)} w_{ij} \left( x_i - x_j \right)  -\kappa_i x_i  + d_i,
\]
where $\kappa_i$ is the degree of stubbornness of node $i$, i.e., the weight that it gives to its own state.
The dynamics of the entire system can be written as:
\begin{equation} \label{ncdyn.eq}
\dot{\xvec} = -(L  + D_{\kappa} D_{S}) \xvec + \dvec,
\end{equation}
where $\dvec$ is a vector of zero-mean white noise processes that affect all nodes, $D_{\kappa}$ is the diagonal matrix of degrees of stubbornness, and $D_{S}$ is a diagonal (0,1) matrix
with its $(i,i)^{th}$ entry equal to 1 if node $i$ is a leader and  0 otherwise.  We note that if $S \neq \emptyset$, then $L + D_{\kappa}D_S$ is positive semi-definite~\cite{RMME09}.

As with noise-free leaders, we define the performance of the system for a given set of leaders $S$ in terms of the total steady-state variance, which is given by~\cite{LFJ14},
\begin{equation}\label{ncvar.eq}
\Rnc(S) = \frac{1}{2} \trace{ (L + D_{\kappa}D_{S})^{-1}}.
\end{equation}
The \emph{noise-corrupted leader selection problem} is to identify the set of at most $k$ leaders that minimizes this variance, i.e.,
\begin{equation}
\begin{array}{ll}
\mathrm{minimize}& \Rnc(S) \\
\mathrm{subject~to}& |S|  \leq k.
\end{array} \label{kleadernc.eq}
\end{equation}

\section{Relationship to Resistance Distance} \label{resistance.sec}

For a graph $G = (V,E,W)$, consider an electrical network with  $V$ the set of nodes and $E$ the set of edges, where each edge $(i,j)$ has resistance $\frac{1}{w_{ij}}$.
 The \emph{resistance distance} between two nodes
$i$ and $j$, denoted $r(i,j)$, is the potential difference between  $i$ and $j$ when a unit current is applied between them.  
Let $L_{j}$ denote the Laplacian matrix 
of $G$ where the row and column of node $j$ has been removed.  It has been shown that~\cite{KR93},
\begin{equation} \label{rL.eq}
r(i,j) = L_{j}^{-1} (i,i),
\end{equation}
i.e., $r(i,j)$ is given by the $(i,i)^{th}$ component of $L_j^{-1}$.

We now show how the performance measures $\Rnf(S)$ and $\Rnc(S)$ can be  expressed  in terms of resistance distances.

\subsection{Noise-Free Leaders}
For a single noise-free leader $v$, it follows directly from (\ref{rL.eq}) that the total steady-state variance is determined
by the resistance distances from all follower nodes to leader node $v$,

\[
\Rnf(\{v \} ) = \frac{1}{2}  \sum_{i \in V \setminus \{v\}} r(i,v).
\]
This relationship can be generalized to multiple noise-free leaders.
In this case,  the resistance distance $r(i,S)$ is the potential difference between follower node
$i$ and the leader set $S$ with unit current.
\begin{proposition} \label{Rr.prop}
The resistance distance $r(i,S)$ from a node $i \in V \setminus S$ to a leader set $S \neq \emptyset$ is related to $\Lf$ as:
\[
r(u,S) = \Lf^{-1}(i,i).
\]
\end{proposition}
\begin{IEEEproof}
Let $B \in \mR^{|E| \times |V|}$ be the incidence matrix of $G$. For each edge $e=(i,j) \in E$, a direction is assigned arbitrarily.
$B(e,i) = 1$ if node $i$ is the tail of  edge $e$, $B(e,i)=-1$ if node $i$ is the head of  edge $e$, and $B(e,i)=0$ otherwise. A resistance $r$ is assigned to each edge $e=(i,j)$ such that $r(e)=\frac{1}{w_{ij}}$. Let  $K \in \mR^{|E| \times |E|}$ be a diagonal matrix with $K(e,e)=r(e)$. It is easy to verify that $B^\top K^{-1} B=L$.
Let $\mathbf{i} \in \mR^{|E|}$ represent the current across all edges, and let $\mathbf{v} \in \mR^{|V|}$ represent the voltages at all vertices. 
By Kirchoff's law, $B^\top \mathbf{i}=\mathbf{c}$,
where $\mathbf{c} \in \mR^{|V|}$ denotes the external currents injected at all vertices, and by Ohm's law, $K\mathbf{i}=B\mathbf{v}$.
It follows that,
\begin{equation} \label{Lv.eq}
L\mathbf{v}=\mathbf{c}.
\end{equation}

Let $\mathbf{v}_j=0$ for all leaders $j \in S$, and thus $\mathbf{v}\tp=[\mathbf{0} ~\mathbf{v}_f\tp]$, where $\mathbf{v}_f$ denotes the voltages for the follower nodes. 
Let $\mathbf{c}_i =1$ for follower $i$ and  ${\mathbf{c}}_k =0$ for followers $k\neq i$. 
 Expanding (\ref{Lv.eq}), we obtain,
\begin{equation}
\left(
\begin{array}{cc}
\Lf & \Llf \\
\Lfl &\Lf
\end{array}
\right)
\left(
\begin{array}{c}
\mathbf{0} \\
\mathbf{v}_{f}
\end{array}
\right)
=
\left(
\begin{array}{c}
\mathbf{c}_l \\
\mathbf{e}_{i}
\end{array}
\right),\nonumber
\end{equation}
where $\mathbf{e}_{i}$ is the canonical basis vector.
Therefore, $\Lf \mathbf{v}_f=\mathbf{e}_i$, 
 Since $\Lf$ is positive definite, and thus, invertible, we have $r(i,S)=\mathbf{v}_i= \Lf^{-1}(i,i)$. 
\end{IEEEproof}

The coherence for a set of noise-free leaders is given in the following theorem, which follows directly from Proposition~\ref{Rr.prop} and (\ref{coherence.eq}).
\begin{theorem}
Let $G$ be a network with  noise-free leader dynamics, and let $S$ be the set of leaders. The coherence of $G$ is:
\[
\Rnf(S) = \frac{1}{2} \sum_{i \in V \setminus S} r(i,S).
\]
\end{theorem}

\begin{figure} 
\centering
\includegraphics[scale=.35]{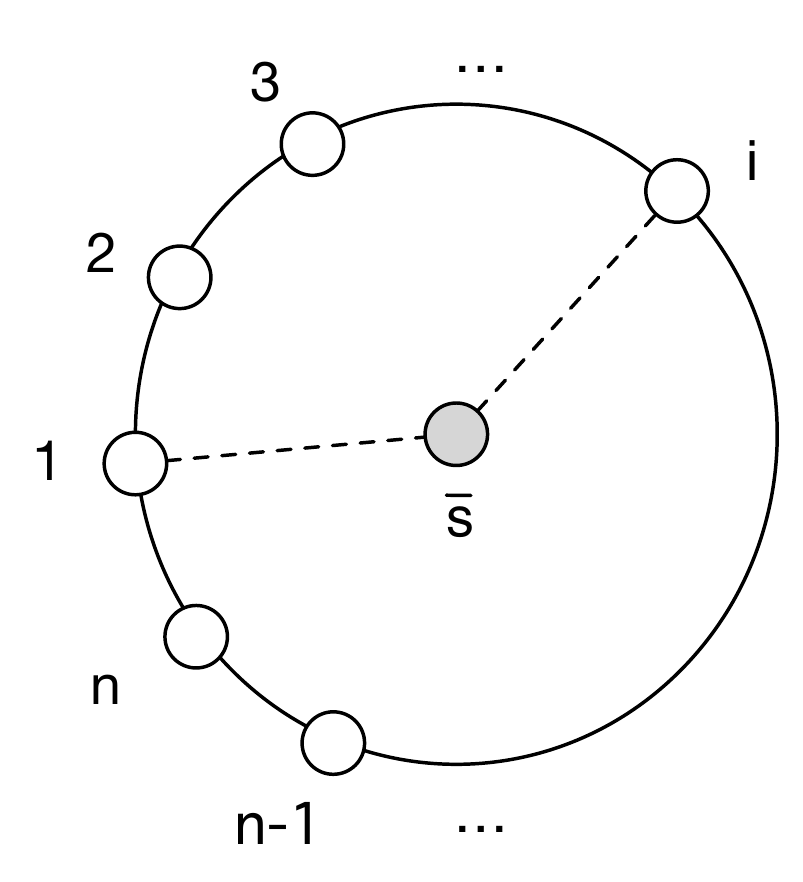}
\caption{Augmented cycle graph with noise-corrupted leaders nodes 1 and $i$.} \label{leadercycle.fig}
\end{figure}

\subsection{Noise-Corrupted Leaders}
For the case of noise-corrupted leaders, we obtain our expression for the coherence by constructing an augmented network.
Let $G = (V,E,W)$ be undirected weighted graph, and let $S \subseteq V$ be a set of noise-corrupted leaders.
We form the augmented graph $\Gbar$ from $G$ by adding a single node $\sbar$ to $G$ and creating an edge from each node $i \in S$ to $\sbar$, with edge weight $\kappa_i$.
An example is shown in Fig.~\ref{leadercycle.fig} for an $n$-node cycle. The noise-corrupted leaders are nodes 1 and $i$.
 We let $\overline{r}(u,v)$ denote the resistance distance between nodes $u$ and $v$ in $\Gbar$

The relationship between resistance distances in $\Gbar$ and the coherence with a set $S$ of noise-corrupted leaders is given in the following theorem.
\begin{theorem} \label{nc.thm}
Let $G=(V,E,W)$ be a network with noise-corrupted leader dynamics, and  let $S$ be the set of leaders.  Let $\Gbar=(\Vbar, \Ebar, \Wbar)$ be the corresponding augmented graph. 
Then, the coherence of $G$ is:
\[
\Rnc(S) =  \frac{1}{2} \sum_{i \in V} \overline{r}(i,\sbar).
\]
\end{theorem}
\begin{IEEEproof}
Let $L$ be the weighted Laplacian of $G$, and let $\overline{L}$ be the weighted Laplacian of $\Gbar$.  We denote by $\overline{L}_{\overline{s}}$  the matrix formed from $\overline{L}$ by removing the row and column corresponding to node $\sbar$.  
We first note that, by the construction of $\Gbar$, $\overline{L}_{\overline{s}} = L + D_S D_{\kappa}$.
By (\ref{rL.eq}), for any node $i \in \Vbar$,
\[
\overline{r}(i, \sbar) = \left(L + D_S D_{\kappa}\right)^{-1}(i,i),  
\]
from which we obtain:
\begin{align*}
\sum_{i \in V} \overline{r}(i,\sbar) &= \trace{(L + D_S D_{\kappa})^{-1}} = 2\Rnc(S),
\end{align*}
where the last equality follows from (\ref{ncvar.eq}).
\end{IEEEproof}

\subsection{Useful Results on Resistance Distance} \label{useful.sec}

We conclude this section by stating some useful results on resistance distance.

\begin{proposition} \label{directS.prop}
Let $S \subseteq V$, and for a  node $i \subseteq V$, let $U_i \subseteq S$ be the set of nodes in $S$ for which there is path from 
$i$ to some $j \in U$ that does not traverse any other element in $S$.  Then $r(i,S) = r(i,U_i)$.
\end{proposition}
This proposition follows directly from the definition of resistance distance.

\begin{lemma}[\cite{KR93} Thm. D] \label{graphDist.lem}
Consider an  undirected connected graph $G = (V,E,W)$, and let $d_{uv}$ denote the graph distance between $u,v \in V$, 
i.e., the sum of the edge weights along the shortest path between $u$ and $v$. Then, $r(u,v) \leq d_{uv}$, with equality if and only if there is single path between $u$ and $v$.
\end{lemma}

\begin{lemma} \label{cutpoint.thm}
Consider a weighted, undirected graph $G=(V,E,W)$, partitioned into two components $A = (V_A, E_A, W_A)$ and $B = (V_B, E_B, W_B)$ that share only a single vertex $\{x\}$.
Let $S \subseteq V_B$.  Then for any $u \in V_A$, 
\[
r(u,S) = r(u,x) + r(x,S)
\]
\end{lemma}
 Lemma~\ref{cutpoint.thm} is a generalization of Lemma~E from~\cite{KR93}.

\begin{lemma} \label{rpRes.lem}
Consider a weighted, undirected path graph, with end vertices $x$ and $y$.  Let $u$ be a vertex on the path.  For any vertices $i,j$ on the path, let $d_{ij}$ denote
their graph distance. Then,
\begin{equation}
r(u,\{x,y\})=d_{ux}-\frac{d_{ux}^2}{d_{xy}}= d_{uy}-\frac{d_{uy}^2}{d_{xy}}.
\end{equation}
\end{lemma}
\begin{IEEEproof}
We start by assigning $0$ voltage to $x$ and $y$. Then, we impose unit external current to $u$, which flows out from $\{x,y\}$.  By definition,  $r(u,\{x,y\}) = \mathbf{v}_u$,
as defined in the proof of Proposition~\ref{Rr.prop}.
It follows that,
\begin{align}
r(u,\{x,y\}) &= \frac{1}{\frac{1}{r(u,x)}+\frac{1}{r(u,y)}} = \frac{1}{\frac{1}{d_{ux}}+\frac{1}{d_{uy}}}   \label{rp2.eq} \\
&= d_{ux}-\frac{d_{ux}^2}{d_{xy}}= d_{uy}-\frac{d_{uy}^2}{d_{xy}}, \label{rp3.eq}
\end{align}
where the second equality follows from Lemma~\ref{graphDist.lem}, and (\ref{rp3.eq}) is obtained from (\ref{rp2.eq}) by applying the equality ${d_{ux}+d_{uy}=d_{xy}}$.
\end{IEEEproof}



\begin{theorem}[\cite{YK13}, Thm. 2.1]  \label{recres.thm}
Let $G'=(V,E',W')$ be the graph formed by adding edge $(i,j)$ to the connected, undirected graph $G=(V,E,W)$, with edge weight $w_{ij}$.
For $p,q \in V$, let $r(p,q)$ denote their resistance distance in $G$, and let $r'(p,q)$ denote their resistance distance in $G'$.  Then,
\[
r'(p,q) = r(p,q) - \frac{w_{ij} \left[ r(p,i) + r(q,j) - r(p,j) - r(q,i)\right]^2}{4\left( 1 + w_{ij} r(i,j)\right)}.
\]
\end{theorem}


\section{Leader Selection Analysis} \label{graphs.sec}

In this section, we use the resistance distance based formulations for coherence to provide closed-form solutions to the leader selection problems for several classes of networks.

We first consider the case of a single leader $v$.
For the noise-free case,
\begin{equation}\label{singleleader.eq}
\Rnf(\{v\}) = \frac{1}{2} \sum_{u \in V \setminus \{v\}} r(u,v).
\end{equation}
The expression (\ref{singleleader.eq}) shows that the optimal single noise-free leader is the node with minimal total resistance distance to all other nodes.
As shown in~\cite{FL16}, this corresponds to the node with  maximal information centrality.

In the noise-corrupted case, 
\begin{align}
\Rnf(\{v\}) &= \frac{1}{2} \sum_{u \in V} \overline{r}(u,\overline{s}) \nonumber \\
 &= \frac{1}{2} \left( \left( \sum_{u \in V \setminus \{v\}} r(u,v) \right) + |V| \kappa_v \right). \label{singleleadernc.eq}
\end{align}
where the last equality follows from Lemma~\ref{cutpoint.thm}.
If all nodes exhibit the same degree of stubbornness, then the optimal noise-free leader and the optimal noise-corrupted leader coincide.
However, if nodes exhibit different degrees of stubbornness, the single best leader may differ for the two dynamics.

We next explore the leader selections problems for $k >1$ leaders.
For the remainder of this section, we restrict our study to networks where all edge weights and all degrees of stubbornness, $\kappa_i$, are equal to 1.

\subsection{$k$ Noise-Free Leaders in a Cycle}

Consider a cycle $n$ nodes, identified by $1, 2, \ldots, n$ in a clockwise direction.
We use the notation $x \prec y$ to mean that node $x$ precedes node $y$ on the ring, clockwise
and $d_{xy}$ denotes the graph distance between nodes $x$ and $y$ where $x \prec y$.
For example, for $n=5$, where $x =4$ and $y=2$, $x \prec y$ and $d_{xy} = 3$.

\begin{theorem} \label{ring.thm}
Let $G$ be an $n$-node cycle with $k$ noise-free leaders, with $n$ written  $n=k\ell + q$, where $\ell$ and $q$ are integers with $0 \leq q < k$.
Let $S = \{s_1, \ldots, s_{k}\}$ be the leaders and let $c$ be a $k$-vector of graph distances between adjacent leaders, i.e., $c_i$ is the distance between leader $s_{i}$ and leader $s_{i+1}$, for $i=1\ldots k-1$, and $c_k$ is the distance between leader $s_k$ and leader $s_1$.
Then:
\begin{enumerate}
\item The coherence of $G$ is 
\[
\Rnf(S) = \textstyle \frac{1}{12}( c\tp c - k).
\]
\item $S$ is an optimal solution to the $k$-leader selection problem if and only if $c \in \cC$, where
\begin{equation*}
\cC = \left\{ d~|~d_i \in \{\ell, \ell+1\},~i=1 \ldots k, \sum_{i=1}^{k} d_i = n \right\}.
\end{equation*}
\end{enumerate}
\end{theorem}

\begin{IEEEproof}
We first find the total resistance  distance to $S$ for all nodes $j$ such that with $s_i \prec j \prec s_{i+1}$,
\begin{align}
\sum_{s_i \prec j \prec s_{i+1}} r(j, S) &= \sum_{s_i \prec j \prec s_{i+1}} r(j, \{s_i, s_{i+1} \}) \label{ringsum1.eq} \\
&= \sum_{\ell=1}^{c_i-1} \textstyle \left( \ell - \frac{\ell^2}{c_i} \right) \label{ringsum2.eq} \\
&= \textstyle \frac{1}{6}\left(c_i^2 - 1\right), \label{ringsum3.eq}
\end{align}
where (\ref{ringsum1.eq}) follows from (\ref{ringsum2.eq}) by Proposition~\ref{directS.prop} and Lemma~\ref{rpRes.lem}.
Applying Proposition~\ref{Rr.prop}, we obtain,
\[
\Rnf(S) = \frac{1}{2} \sum_{i=1}^k \frac{1}{6}( c_i^2 -1) = \frac{1}{12} (c\tp c - k).
\]
With this, we can express  problem (\ref{kleader.eq}) as an integer quadratic program:
\begin{align}
\text{minimize}~&\textstyle \frac{1}{12} \left( c\tp c - k \right) \label{mathobj.eq} \\
\text{subject to}~&\textbf{1}\tp c = n \label{sumcon.eq} \\
~& c_i \in \{1, \ldots, n-1\}. \label{intcon.eq}
\end{align}
If $k$ divides $n$, then it is straightforward to verify that ${c^{\star} = \frac{n}{k} \textbf{1}}$ is a solution to the above problem.
In this case, $\ell = \frac{n}{k}$.


For $n=kl + q$, with $q > 0$,
assume that $c$ is a solution to (\ref{kleader.eq}) but $c \notin \cC$.
Then there exists some component $c_i$ such that $c_i = \ell+1 + x$ for some integer $x > 0$ and some component $c_j$ such that $c_j = \ell - y$
for some integer $y > 0$.
Let $c'$ be such that $c'_i = c_i - 1$ and $c'_j = c_j + 1$ and $c'_{\ell} = c_{\ell}$ for all $ \ell \neq i, \ell \neq j$.  
Clearly, $\frac{1}{12} \left( c\tp c - k \right) > \frac{1}{12} \left( (c')\tp c' - k \right)$, which contradicts our assumption
that $c$ is a solution to (\ref{kleader.eq}).
\end{IEEEproof}

\subsection{$k$ Noise-Free Leaders in a Path}

Consider a path graph with $n$ nodes, identified by $1, 2, \ldots, n$.  
Let $d_{uv}$ denote the graph distance between nodes $u$ and $v$.

\begin{theorem} \label{path.thm}
Let $G$ be a path graph with $n$ nodes, and let $S = \{s_1, \ldots s_k\}$ be a set of $k$ noise-free leaders.  
Let $c$ be a {$(k+1)$-vector}, where $c_1 = (s_1 - 1)$ and $c_{k+1} = n- s_k$.
Let $c_i = s_{i+1} - s_i$, for $i=2, \ldots, k$. Then, 
\begin{enumerate}
\item The coherence of $G$ is:
\begin{equation}
\Rnf(S) = \frac{1}{4}\left(c_1^2 + c_{k+1}^2 + c_1 + c_{k}\right) + \frac{1}{12} \sum_{i=2}^k (c_i^2 -1). \label{pathcoh.eq}
\end{equation}
\item Let $n$ be such that, for the optimal leader configuration, it holds that $c_1 + c_{k+1} = a$,
where $2$ divides $a$ and $b=(n-1)-a$, where $(k-1)$ divides $b$.  Then, the optimal solution to the $k$-leader selection problem
is:
\begin{align}
c_1 &= c_{k+1} = \text{round}\left( \frac{2(n-1) - 3(k-1)}{6(k-1) + 4} \right)  \label{cpath1.eq} \\
c_i &= \frac{1}{k-1}\left( (n-1) - 2c_1 \right), ~~i=2 \ldots k. \label{cpath2.eq}
\end{align}
\end{enumerate}
\end{theorem}

\begin{IEEEproof}
We first find the total resistance distance to $S$ for all nodes $u$ with $1 \leq u < s_1$:
\begin{equation}
\sum_{u=1}^{s_1 - 1} r(u,S) = \sum_{u=1}^{c_i} u = \frac{c_i(c_i+1)}{2}, \label{endcom1.eq}
\end{equation}
where the first equality follows from Proposition~\ref{directS.prop} and Lemma~\ref{graphDist.lem}.
Similarly, the total resistance distance to $S$  for all nodes $v > s_{k}$ is:
\begin{equation}
\sum_{v=s_{k+1}}^n r(v,S) = \frac{c_{k+1} (c_{k+1} + 1)}{2}. \label{endcom2.eq}
\end{equation}

The total resistance distances to $S$ for all nodes $u$ between $s_i$ and $s_{i+1}$ can be obtained in a similar fashion to (\ref{ringsum1.eq}) - (\ref{ringsum3.eq}),
\begin{equation} \label{midcom1.eq}
\sum_{u=s_i + 1}^{s_{i+1} - 1} r(u,S) = \frac{1}{6}(c_{i+1}^2 - 1).
\end{equation}

Combining (\ref{endcom1.eq}), (\ref{endcom2.eq}), and (\ref{midcom1.eq}) with Proposition~\ref{Rr.prop}, we obtain, 
\begin{small}
\begin{align*}
\Rnf(S) &= \frac{1}{2} \sum_{u \in V \setminus S} r(u,S)\\
&= \frac{1}{4}(c_1^2+  c_{K+1}^2+ c_1 + c_{k+1}) + \frac{1}{12} \sum_{i=2}^k (c_i^2 - 1).
\end{align*}
\end{small}

To find the optimal leader locations, we must solve the optimization problem,
\begin{equation} \label{lineopt.eq}
\begin{array}{ll}
\text{minimize}~~ &  c\tp P c + r\tp d  - \frac{k-1}{12}  \\
\text{subject to}~&\textbf{1}\tp c = n-1 \\
& c_i \in \{1, \ldots, n-1\},
\end{array}
\end{equation}
where $P$ is the $(k+1)\times(k+1)$ diagonal matrix with diagonal components $[ \frac{1}{4}~~\frac{1}{12}~ \ldots ~\frac{1}{12} ~~\frac{1}{4}]$,
and $r$ is a $(k+1)$-vector with $r_1 = r_{k+1} = \frac{1}{4}$, and all other entries equal to 0.
Let $\cstar$ be a solution to (\ref{lineopt.eq}).
Using a similar argument to that in the proof of Theorem~\ref{ring.thm}, we can conclude 
$\cstar_1 = \cstar_{k+1} = a/2$ for some even integer $a$, i.e., 
that leaders $s_1$ and $s_k$ should each be the same distance from their respective ends of the path.
Similarly, $\cstar_i = b/(k-1)$ for $i=2 \ldots k$, i.e.,  the leaders $s_2, \ldots, s_{k-1}$ should be equidistant. 

Let $q = a/2$, so that the optimal leader placement has $\cstar_1 = \cstar_{k+1} = q$ and $\cstar_i = \frac{1}{k-1}((n-1) - 2q)$, for $i=2 \ldots k$.
Then, we can reframe (\ref{lineopt.eq}) as $\text{minimize}_{q \in \{1, \ldots, n\}} C(q)$, where, 
\begin{align*}
C(q) &= \frac{1}{2}(q^2 + q) + \frac{(k-1)}{12}\left( \left( \frac{(n-1)-2q}{k-1} \right)^2- 1\right).
\end{align*}
Relaxing the integer constraint, the value of $q$ that minimizes $C(q)$ is,
\[
q^{\star} = \frac{2(n-1) - 3(k-1)}{6(k-1) + 4}.
\]
Since $C(q)$ is quadratic, the optimal integer value for $q$ is $\text{round}(q)$.
The values for $c_i$, $i =1 \ldots k+1$, in (\ref{cpath1.eq}) - (\ref{cpath2.eq}) follow from the definition of $q$ above.
\end{IEEEproof} 

While the restriction that $c_{1} = c_{k+1}$ and $c_i = c_{i+1}$, $i=2 \ldots k$ does not hold for all network sizes,
it can be shown experimentally to hold for many.  An example is a path graph with $n=40$ and $k=3$, where $c_1 = c_{4} = 3$
and $c_1 = c_2 = c_3 = 11$.

\subsection{Two Noise-Free Leaders in Trees}

We next consider the 2-leader selection problem in rooted, undirected $\M$-ary trees. An \emph{$\M$-ary tree} is a rooted tree where each node has at most $\M$ children. A \emph{perfect $\M$-ary tree} is an $\M$-ary tree in which all non-leaf nodes have exactly $\M$ children and all leaves are in the same level.
Let $r$ denote the root node of the tree, and let $h$ denote its height.  
We number the levels of the tree starting with the root, as $0, 1, 2, \ldots, h.$ The root of the tree is at level 0, and the leaves of a perfect $M$-ary tree of height $h$ are  at level $h$.  We use $\lev{x}$ to denote the level of a node.

%
We begin with the following lemma, which gives general guidance for the optimal location of two leader nodes.
\begin{lemma} \label{leaderSpan.prop}
Consider a perfect $\M$-ary tree  $T=(V,E)$.  Let $x, y \in V, x \neq y$ be such that their lowest common ancestor is a node of level $\ell > 0$.
Then, there exists $y, z \in V, y \neq z$, with lowest common ancestor $r$ such that  $\Rnf(\{x,y\}) > \Rnf(\{y,z\})$.
\end{lemma}
The proof of this lemma is given in Appendix \ref{leaderSpan.prop}.
Lemma~\ref{leaderSpan.prop} tells us that the optimal 2-leader set  will not have two nodes in the same subtree of a child of $r$.

We denote these two leaders by $x$ and $y$, and assume there lowest common ancestor is $r$.  Without loss of generality, we assume $\lev{x} \leq \lev{y}$.
We denote the graph distances between $x$ and $y$, $x$ and $r$, and $y$ and $r$ by $d_{xy}$, $d_{xr}$ and $d_{yr}$, respectively.
To study the coherence of this system, we decompose the tree into three subgraphs, (1) the subtree of $T$ rooted at $y$,  denoted  $T_y = (V_y, E_y)$,
(2) the subtree of $T$ rooted at $x$, excluding those nodes in $T_y$,  denoted by $T_x = (V_x, E_x)$, and (3) the induced subgraph of $T$ consisting of nodes $V - (V_x \cup V_y) \cup \{x,y\}$, which is denoted by $G_{xy} = (V_{xy}, E_{xy})$.  
Note that by Proposition~\ref{directS.prop}, for $u \in V_x$, it holds that $r(u,S) = r(u,x)$.
Similarly, for $u \in V_y$, we have $r(u,S) = r(u,y)$.
We can therefore decompose $\Rnf(S)$ as,
\begin{align} 
&\Rnf(\{x,y\})= \nonumber \\
&~ \frac{1}{2} \left(\sum_{u \in V_x} r(u,x) + \sum_{u \in V_y} r(u,y) + \sum_{u \in V_{xy}} r(u, \{x,y\}) \right)\label{Rtree1.eq} \\
~~~&=  \frac{1}{2}\left(\sum_{u \in V_x} d_{ux} + \sum_{u \in V_y} d_{uy} + \sum_{u \in V_{xy}} r(u, \{x,y\})\right) \label{RofTree.eq},
\end{align}
where (\ref{RofTree.eq}) follows from (\ref{Rtree1.eq}) by Lemma~\ref{graphDist.lem}.

With this decomposition, we can apply the building blocks described in Section~\ref{useful.sec} to identify the optimal 2 noise-free leaders in $M$-ary trees for various values of $M$.
We begin with $M=2$.
\begin{theorem} \label{btree.thm}
For the noise-free 2-leader selection problem in a perfect binary tree with height $h\geq 4$, the optimal leaders are such that  $d_{xy}=4$ and $d_{xr}=2$, and the resulting coherence is: 
\begin{equation}
\Rnf(S)= \frac{(n+1)}{2}\left(\log_{2}{(n+1)}-\frac{25}{8}\right)+\frac{7}{2}\,. \label{Rbtree.eq}
\end{equation}
\end{theorem}
The proof of Theorem~\ref{btree.thm} is given in Appendix~\ref{btreeproofs.sec}.


It is interesting to note that the optimal leader locations are independent of the height of the tree.  This independence of the height also holds for $M >2$, as shown in the following theorems.
Proofs are given in Appendix~\ref{btreeproofs.sec}.
\begin{theorem}
\label{3Theo}
For the noise-free 2-leader selection problem in a perfect ternary tree $T^{(3)}$ with  height $h\geq 4$, the optimal leaders are such that $d_{xy}=2$ and $d_{xr}=1$, and the resulting coherence is:
\begin{equation}
\Rnf(S)=\frac{2n+1}{4}\left(\log_{3}(2n+1)-2\right)+1\,.
\label{Rttree.eq}
\end{equation}
\end{theorem}

\begin{theorem}
\label{KTheo}
For the noise-free 2-leader selection problem a perfect $\M$-ary tree $T$, with $M \geq 4$ and $h \geq 4$, the optimal leaders are such that $d_{xy}=1$ and $d_{xr}=0$, and the resulting coherence is:
\begin{small}
\begin{align}
&\Rnf(S)=\frac{1}{2}\left(n+\frac{1}{M-1}\right)\log_{M}(nM-n+1)\nonumber\\
&\quad -\frac{n(M^2+M-1)}{2M(M-1)}+\frac{1}{2M}\,.
\label{Rktree.eq}
\end{align}
\end{small}
\end{theorem}

\subsection{Two Noise-Corrupted Leaders in a Cycle Graphs}

Consider an $n$-node cycle with nodes  labeled $\{1,2,\ldots,n\}$. 
We use Theorem~\ref{recres.thm} to determine the coherence of the graph as a function of the graph distance between nodes 1 and $i$.

\begin{theorem} \label{ncring.thm}
In an $n$-node cycle with two noise-corrupted leaders,  where $n$ is even, the coherence is minimized with the leaders are at distance $n/2$ apart, and the resulting coherence is:
\begin{equation}
\Rnc(S) = \frac{n^3 + 16n^2 + 44n - 16}{24 (n + 8)}. \label{Rnc.eq}
\end{equation}
\end{theorem}

\begin{IEEEproof}
Without loss of generality, we assume node 1 and node $i$ are noise-corrupted leaders.
Let the graph $\Gbar$ be the augmented graph shown in Fig.~\ref{leadercycle.fig}, omitting edge $(i,\overline{s})$.
By Lemma~\ref{rpRes.lem}, for  arbitrary nodes $u,v \in \{1,2,...,n\}$, their resistance distance in $\Gbar$ is:
\[
\overline{r}(u,v) = \frac{| u-v | ( n - |u-v|)}{n}, 
\]
and the resistance distance from a node $u \in \{1,2,\ldots,n\}$ is:
\[
\overline{r}(u,s) = r(u,1) + 1 = \frac{(u-1)(n- (u-1))}{n} + 1.
\]

Let $G'$ be the graph formed from $\overline{G}$ by the addition of edge $(i,\overline{s})$.  Then, for a node $u \in \{1,2,...,n\}$, the resistance distance from $u$ to $\sbar$ in $G'$ is: 
\begin{align*}
&r'(u,\sbar) = \rbar (u,\sbar) - \frac{\left[ \rbar(u,i) + r(\sbar,\sbar) - \rbar(u,\sbar) - \rbar(i,\sbar) \right]^2}{4\left( 1 + \rbar(i,\sbar)\right)}  \\
&= \frac{(u-1)(n-(u-1)}{n} + 1 - \\
& \frac{\left[\frac{|u-i| (n- |u-i|)}{n}  - \frac{(u-1)(n-(u-1))}{n}  - \frac{((i-1)(n-(i-1))}{n} -2 \right]^2}{4\left( 1 + \left(\frac{(i-1)(n-(i-1))}{n} + 1\right)\right)}.
\end{align*}

By Theorem~\ref{nc.thm}, summing over all nodes $u$, we obtain:
\begin{align}
&\Rnc(S) = \frac{1}{2} \sum_{u=1}^n r'(u,\sbar)  = \frac{1}{12}(n^2 + 6n - 1) \nonumber \\ 
&-  \frac{1}{12n\left( 2 + \frac{(i-1)(n-(i-1))}{n} \right)} \Big[2i^4 - 4i^3(n+2) \nonumber \\
&+ i^3(2n^2 + 6n + 11)+ i(2n^2 + n - 6)  + 2n^2 - 3n + 1 \Big]. \label{Rnci.eq}
\end{align}
We note that this function is continuous over the interval $[1,n].$

We then take the derivative with respect to $i$:
\begin{align*}
&\frac{\partial }{\partial i} \Rnc(S) = \\
&\frac{1}{(6 (-i^2 + i (n + 2) + n - 1)^2)}\Big[ 2 i^5 - 5 i^4 (n + 2) )\\
&+ 4 i^3 (n^2 + 3n + 5) - i^2 (n^3 + 6n + 20) \\
&- 2 i (n^3 + 3n^2 + n - 5) + n^2 + n- 2 \Big].
\end{align*}
The derivative has five roots.
Of these, only $i = (n + 2)/2$  lies in the interval $[1,n]$. Further, it is a minima of $\partial / \partial i \Rnc(S)$.
For  even $n$, we substitute this value of $i$ into (\ref{Rnci.eq}) to obtain (\ref{Rnc.eq}).
\end{IEEEproof}

\section{Comparison to Coherence in Leader-Free Networks} \label{compare.sec}

Network coherence has also been studied in graphs without leaders. 
In this setting, every node behaves as a follower, using the dynamics in (\ref{fdyn.eq}).
Coherence is measured as the total steady-state variance of the deviation from the average of all node states,
\[
V = \lim_{t \rightarrow \infty} \sum_{i=1}^n \expec{x_i(t) - \frac{1}{n} \sum_{j=1}^n x_j(t)}^2.
\]
It has been shown that, for a network with a single noise-free leader, i.e.,$|S| = 1$~\cite{PB10},
\[
\Rnf(S) \geq V.
\]
In some sense, this means that adding a single leader increases the disorder of the network.

In a leader-free cycle graph, it has been shown that the coherence $V$ scales as $O(n^2)$~~\cite{BJMP12}.
In a cycle with $k$ noise-free leaders, where the leaders are located optimally, by Theorem~\ref{ring.thm}, 
\[
\Rnf(S) = \frac{1}{12}\left( \left( \frac{n}{k}\right)^2 \textbf{1}\tp \textbf{1} - k\right) = \frac{n^2}{12k} - \frac{k}{12}.
\]
Thus for a fixed leader set size $k$, The coherence $\Rnf(S)$ also scales as $O(n^2)$.
Similarly, for the optimal two noise-corrupted leaders in a cycle, $\Rnc(S)$ scales as $O(n^2)$.
This shows that, in the limit of large $n$, in cycle networks, the disorder of the network is similar for leader-free and leader-follower consensus networks.

\section{Numerical Example} \label{treegrowing.sec}

\begin{algorithm}[t]
\caption{Algorithm to add nodes to a perfect binary tree of height $h$, while maintaining optimality of the 2 leaders.}\label{fillTree.alg}
\begin{algorithmic}\small
\Require{$T_h$, with optimal 2 leaders $x$ and $y$}\\
\While{there is a new node $u$ to add}
\If{last level of left or right subtree of  $x$ is not filled}
	\State{Add node $u$ to level $h+1$ of subtree of $x$ with fewer leaves,}
	\State{breaking ties arbitrarily.}
\ElsIf{last level of left or right subtree of  $y$ is not filled}
	\State{Add node $u$ to level $h+1$ of subtree of $y$ with fewer leaves,}	
	\State{breaking ties arbitrarily.}
\Else
	\State{Add node $u$ as leaf on level $h+1$, in any  remaining location.}	
\EndIf
\If{level $h+1$ is filled}
	\State{$h \gets h+1$}
\EndIf
\EndWhile
\end{algorithmic}
\end{algorithm}

\begin{figure}
\centering
\includegraphics[scale=.4]{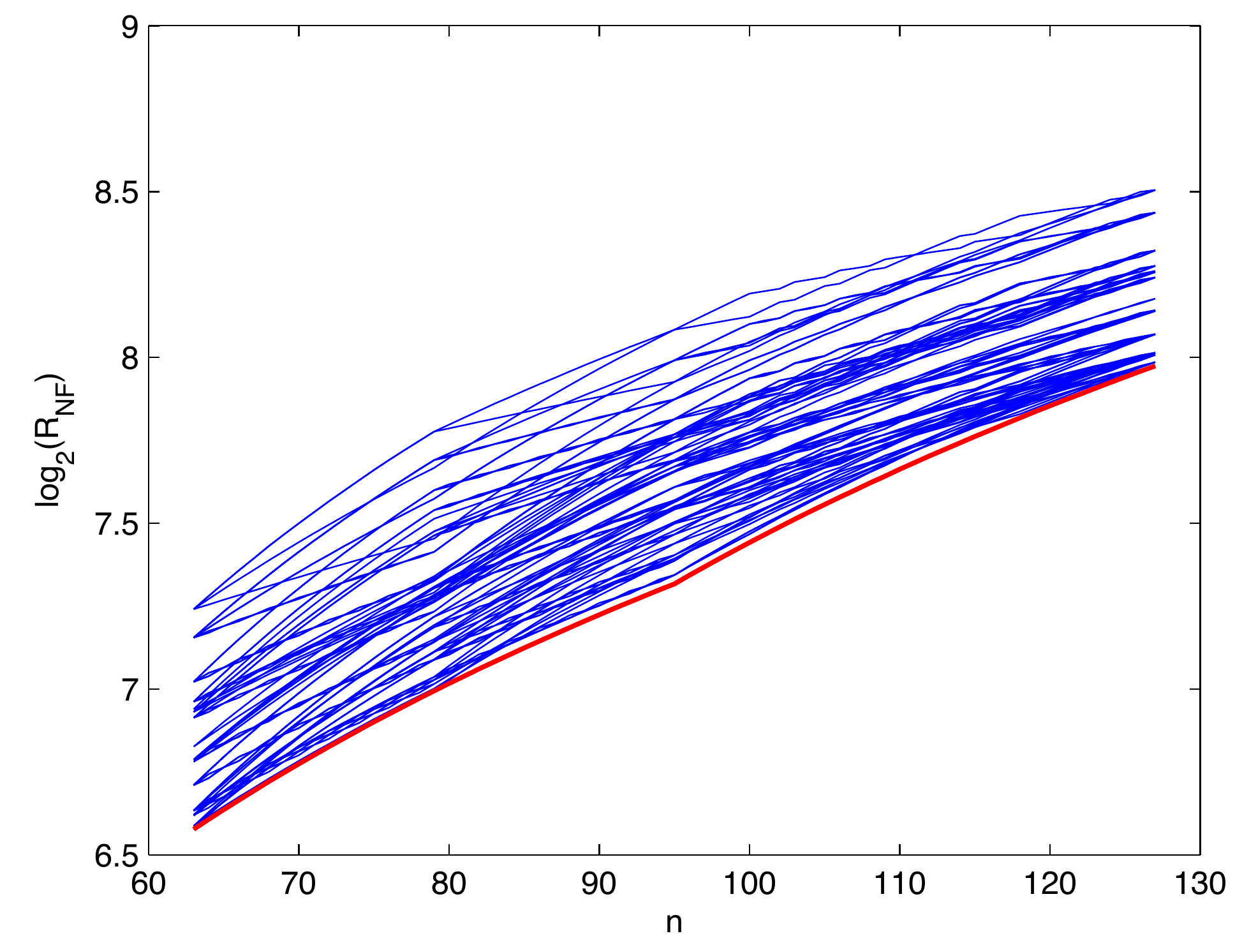}
\caption{Coherence for two-noise free leaders, as a perfect binary tree  of height 5 is grown into a perfect binary tree of height 6 using Algorithm~\ref{fillTree.alg}.} \label{logRnfSim.fig}
\end{figure}

Theorem~\ref{btree.thm} applies to the noise-free leader selection problem in a perfect binary tree.  We now present an algorithm that can be used to grow the tree by adding nodes in a way that does not change the location of the optimal two noise-free leaders.
Pseudocode for this tree-growing process is given in Algorithm~\ref{fillTree.alg}. The algorithm is initialized with a perfect binary tree of height $h \geq 4$, with the optimal leader set $\{ \hat{x}, \hat{y}\}$, with $d_{\hat{x}r} = 2$ and $d_{\hat{y}r} = 2$.  In each iteration, a node is added in a location dictated by the algorithm.

The analysis of this algorithm remains an open question. However, example executions, the algorithm is able to grow a tree from height $h$ to height $h+1$ without impacting the optimality of the leader nodes $x$ and $y$.
In Fig.~\ref{logRnfSim.fig}, we show such and execution. The algorithm is initialized with a perfect binary tree of of height $h = 5$, with 63 nodes. Nodes are added according to the algorithm, until the tree is a perfect binary tree of height $h=6$, with 127 nodes.
The figure shows the coherence for every pair of leader nodes such that $d_{xr} \leq 3$ and $d_{yr} \leq 3$, in log scale.  The coherence for the  leader set $\{ \hat{x}, \hat{y}\}$ is shown in red, while the coherence for each other leader set is shown in blue.
As the figure shows, the coherence  for $\{\hat{x}, \hat{y}\}$ is the smallest throughout the execution of the algorithm.

\section{Conclusion} \label{conclusion.sec}
We have investigated the performance of leader-follower consensus networks under two types of leader dynamics, noise-free leaders and noise-corrupted leaders.
For both leader dynamics, we have developed a characterization of the system performance in terms of  resistance distances in electrical networks.
With this characterization, we have derived closed-form expressions for network coherence in terms of the leader locations. We have also identified the optimal leader locations
in several special classes of networks. 

In future work, we plan to extend our analysis to study coherence in general leader-follower networks.  We also plan to develop a similar mathematical framework to study coherence in second-order systems.


%

\appendices

\begin{figure}
\begin{center}
\includegraphics[width=.9 \linewidth]{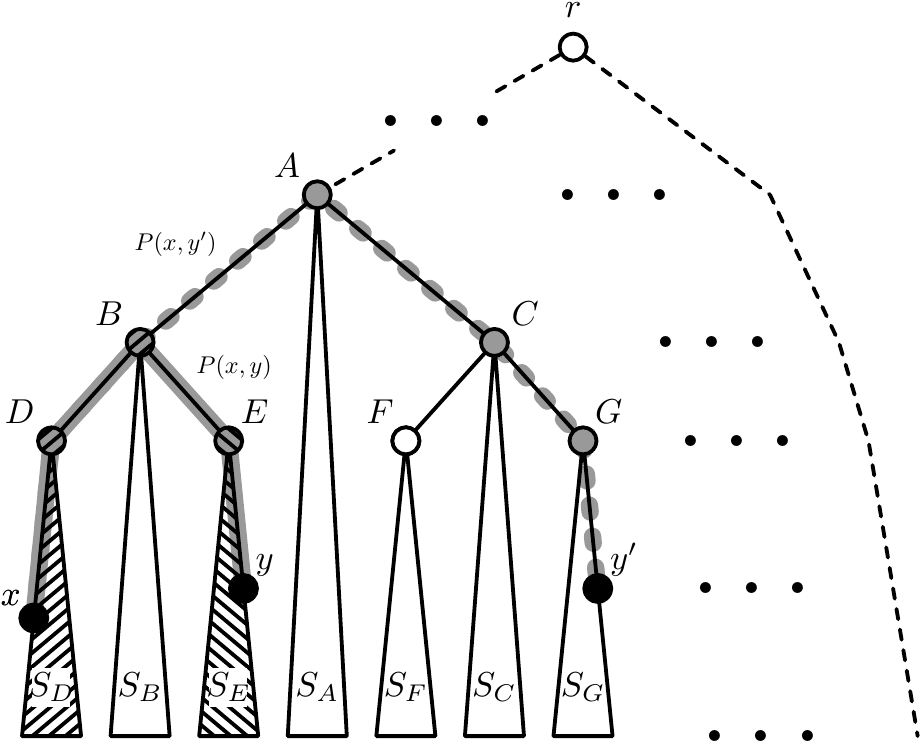}
\caption{Arrangement of nodes in perfect binary tree with two possible leader sets ~$\{x,y\}$ and $\{x,y'\}$.}\label{oneside}
\end{center}
\end{figure}

\section{Proof of Lemma~\ref{leaderSpan.prop}}
\begin{IEEEproof}
Let $x$ and $y$ be the optimal two noise-free leaders in a perfect binary tree $T$.
Without loss of generality, let ${\text{lev}(x) \leq \text{lev}(y)}$.
Assume, for contradiction, that the lowest common ancestor of $x$ and $y$ is a node $B$ that is a descendant of the root $t$.
 Let $A$ be the parent of $B$, and let $D$ and $E$ be children of $B$. 
Let $x$ be a member of the node set consisting of node $B$ and the nodes in the subtree rooted at $D$.  We denote this node set by $S_D$.
Let $y$ be a node in the subtree rooted at $E$. We denote this node set by $S_E$. Let $S_B$ denote the node set in the subtree rooted at $B$, excluding $B$ and the nodes in $S_D$ and $S_E$.
This arrangement is shown in Fig.~\ref{oneside}.

Let $C$ be another child of $A$, as shown in the figure, and let $F$ and $G$ be children of $C$, with the node sets of the the trees rooted at $F$ and $G$ denoted by $S_F$ and $S_G$, respectively.
Let $S_C$ denote the node set of the subtree rooted at $C$, excluding $C$ and the nodes in $S_F$ and $S_G$.
We will prove that, for a node $y'$ in the subtree rooted at $G$ that is in the same location with respect to $G$ that $y$ is with respect to to $E$, 
$\Rnf{\{x,y\}} > \Rnf{\{x, y'\}}$.

Let $P(x,y)$  denote the vertices along the path between $x$ and $y$, and let $P(x,y')$ denote the vertices along the path between $x$ and $y'$.
Consider a pair of vertices $u\in S_D$ and $k\in S_F$, where $k$ is at the same location relative to $F$ (in the subtree rooted at $F$) that $u$ is relative to $D$ (in the subtree rooted at $D$).
Denote the vertex on $P(x,y)$ that is the nearest to $u$ by $p$. 
We find the sum of the resistance distances of $u$ and $k$ to the respective leader sets $\{x,y\}$ and $\{x,y'\}$.
By Lemmas~\ref{graphDist.lem} and \ref{rpRes.lem},
\begin{eqnarray}
R_{u,k} &=& r(u,\{x,y\})+r(k,\{x,y\})\nonumber\\
             &=& d_{u p}+r(p,\{x,y\})+d_{kB}+r(B,\{x,y\})\,, \nonumber \\
R^{\prime}_{u,k} &=& r(u,\{x,y'\})+r(k,\{x,y'\})\nonumber\\
             &=& d_{up}+r(p,\{x,y'\})+d_{kC}+r(C,\{x,y'\})\,.\nonumber
\end{eqnarray}
Noting that $d_{By}=d_{Cy'}$ and $d_{xy'}=d_{xy}+2$, and applying Lemma~\ref{rpRes.lem}, we have:
\begin{align}
R_{u,k}-R^{\prime}_{u,k}=&(d_{up}-d_{up})+(d_{kB}-d_{kC})\nonumber\\
&+\left(r(p,\{x,y\}-r(p, \{x,y^\prime\})\right)\nonumber\\
&+\left(r(B,\{x,y\})-r(C,\{x,y'\})\right)\nonumber\\
=&0+2+\left(d_{px}-\frac{d^2_{px}}{d_{xy}}-(d_{px}-\frac{d^2_{px}}{d_{xy'}})\right)\nonumber\\
&+\left(d_{By}-\frac{d^2_{By}}{d_{xy}}-(d_{Cy'}-\frac{d^2_{Cy'}}{d_{xy'}})\right)\nonumber\\
&=2-\frac{2(d_{xp}^2+d_{By}^2)}{d_{xy}(d_{xy}+2)}\,.\nonumber
\end{align}
Further, since $d_{xp}\leq d_{xB}$ and $d_{xB}+d_{By}=d_{xy}$, we obtain:
\begin{equation}
R_{u,k}-R^{\prime}_{u,k}\geq \frac{4(d_{xB}d_{By}+d_{xy})}{d_{xy}(d_{xy}+2)}>0\,.\nonumber
\end{equation}

Next, consider a pair of vertices $v\in S_E$ and $\ell \in S_G$, where $\ell$ is at the same location relative to $G$ that $v$ is relative to $E$.  
Denote the vertex on $P(x,y)$ that is the nearest to $v$ by $q$, and denote the vertex on $P(x,y')$ that is nearest to $\ell$ by $m$
The sum of the resistance distances from $v$ and $\ell$ to the respective leader sets are (again, by Lemmas~\ref{graphDist.lem} and ~\ref{rpRes.lem}),
\begin{eqnarray}
R_{v,\ell} &=& r(v,\{x,y\})+r(\ell,\{x,y\})\nonumber\\
             &=& d_{vq}+r(q,\{x,y\})+d_{\ell B}+r(B,\{x,y\}) \nonumber \\
    R^{\prime}_{v,\ell} &=& r(v,\{x,y'\})+r(\ell,\{x,y'\})\nonumber\\
             &=& d_{v B}+r(B,\{x,y'\})+d_{\ell m}+r(m,\{x,y'\})\,.\nonumber
\end{eqnarray}


We note that $d_{vq} = d_{\ell m}$, $d_{qy}=d_{my'}$ and $d_{xy'}=d_{xy} + 2$.  Then,
\begin{align}
R_{v,\ell}-R^{\prime}_{v,\ell}=&(d_{vq}-d_{lm})+(d_{lB}-d_{vB})\nonumber\\
&+r(q,\{x,y\})-r(m,\{x,y'\})\nonumber\\
&+r(B,\{x,y\})-r(B,\{x,y'\})\nonumber\\
=&0+2+\left(d_{qy}-\frac{d^2_{qy}}{d_{xy}}-(d_{my'}-\frac{d^2_{my'}}{d_{xy'}})\right)\nonumber\\
&+\left(d_{xB}-\frac{d^2_{xB}}{d_{xy}}-(d_{xB}-\frac{d^2_{xB}}{d_{xy'}})\right)\nonumber\\
=&2-\frac{2(d_{xB}^2+d_{qy}^2)}{d_{xy}(d_{xy}+2)}\,.\nonumber
&\end{align}
Since  $d_{qy}\leq d_{By}$ and $d_{xB}+d_{By}=d_{xy}$, similar to $R_{u,k}-R^{\prime}_{u,k}$, we obtain $R_{v,\ell}-R^{\prime}_{v,\ell}>0$.

For a pair of vertices $i\in S_B\cup \{B\}$ and $j\in S_C\cup \{C\}$, where $\ell$ is at the same location relative to $G$ that $v$ is relative to $E$, define
\begin{align*}
R_{i,j} &= r(i,\{x,y\}) + r(j,\{x,y\})\\
&= d_{iB} + r(B,\{x,y\}) + d_{jB} + r(B,\{x,y\}) \\
R^{\prime}_{i,j} &= r(i,\{x,y'\}) + r(j,\{x,y'\})\\
&= d_{iB} + r(B,\{x,y'\}) + d_{jC} + r(C,\{x,y'\}).
\end{align*}
Since $d_{Cy'}=d_{By}$,
\begin{align*}
R_{i,j}-R^{\prime}_{i,j}=&(d_{iB}-d_{iB})+(d_{jB}-d_{jC})\\
&+r(B,\{x,y\})-r(B,\{x,y'\})\\
&+r(B,\{x,y\})-r(C,\{x,y'\})\\
=&0+2+\left(d_{Bx}-\frac{d^2_{Bx}}{d_{xy}}-(d_{Bx}-\frac{d^2_{Bx}}{d_{xy'}})\right)\\
&+\left(d_{By}-\frac{d^2_{By}}{d_{xy}}-(d_{Cy'}-\frac{d^2_{Cy'}}{d_{xy'}})\right)\\
=&2-\frac{2(d^2_{xB}+d^2_{By})}{d_{xy}(d_{xy}+2)}\,. 
\end{align*}
Since $d_{xB}+d_{By}=d_{xy}$, it follows that
$$R_{i,j}-R^{\prime}_{i,j}=\frac{4(d_{xB}d_{By}+d_{xy})}{d_{xy}(d_{xy}+2)}>0\,.$$ 

Finally, we consider a vertex $t$ that is not in the subtree rooted at $B$ nor the subtree rooted at $C$. 
In this case,
\begin{align*}
R_t &= r(t,\{x,y\}) = d_{tB}+r(B,\{x,y\})\\
R^{\prime}_t &= r(t,\{x,y'\}) = d_{tA}+r(A,\{x,y'\}).
\end{align*}
It follows that
\begin{equation}
R_t-R^{\prime}_t=2-\frac{d_{xB}^2}{d_{xy}}+\frac{(d_{xB}+1)^2}{d_{xy}+2}\,.\nonumber
\end{equation}
Recall that $d_{xy}=d_{xB}+d_{By}$. Thus, $R_t-R^{\prime}_t>0$.

Since $|S_D|=|S_E|=|S_F|=|S_G|$ and $|S_B\cup \{B\}|={|S_C\cup \{C\}|}$, by grouping vertices into pairs, we have shown that,
\[
\frac{1}{2} \sum_{i=1}^n r(i,\{x,y\}) > \frac{1}{2} \sum_{i=1}^n r(i, \{x,y'\}).
\]
This contradicts our assumption that $\{x,y\}$ is the optimal leader set.

\end{IEEEproof}


\section{Proof of Theorems~\ref{btree.thm}, \ref{3Theo}, and \ref{KTheo}} \label{btreeproofs.sec}
We first define a quantity $\Omega(S)$ as:
\begin{align*}
\Omega(S)=\sum_{i\in V\backslash S}r(i,S)=2\Rnf(S)\,.
\end{align*}
and note that a set $S$ that is a minimizer of $\Omega(\cdot)$ is also a minimizer of $\Rnf(\cdot)$.

We next present a lemma that gives $\Omega(\cdot)$ of a perfect $M$-ary tree with two noise free leaders .
\begin{lemma} \label{Mxy.lem}
Let $T$ be a perfect $\M$-ary tree with height $h$. Let $x$ and $y$ be its two noise-free leaders, and assume that the lowest common ancestor
of $x$ and $y$ is the root of $T$. Then,
\begin{small}
\begin{align}
&\Omega(\{x,y\}) = \frac{\M^{h+1}+1}{\M-1}\left(d_{xr}-\frac{d^2_{xr}}{d_{xy}}\right)\nonumber\\
&+M^{h+1}\left(\frac{2}{(M-1)^2}+\frac{M+1}{(M-1)^3dy}\right)(M^{d_{xr}-d_{xy}}+M^{-d_{xr}})\nonumber\\
&+M^{h+1}\left(\frac{h}{M-1}-\frac{3}{(M-1)^2}-\frac{2(M+1)}{(M-1)^3d_{xy}}\right)\nonumber\\
&+\frac{d_{xy}}{M-1}+\frac{M}{(M-1)^2}. \label{rgkary}
\end{align}
\end{small}
\end{lemma}

\begin{IEEEproof}
Recall that in (\ref{RofTree.eq}), we decomposed the coherence into three terms: the coherence in the subtree rooted at $x$, the coherence in the subtree rooted at $y$, and the coherence at the remaining nodes. We can also devide $\Omega$ into three part as
\begin{align*} 
&\Omega(\{x,y\})=  \sum_{u \in V_x} d_{ux} + \sum_{u \in V_y} d_{uy} + \sum_{u \in V_{xy}} r(u, \{x,y\}) .
\end{align*}
We let $T_y$ denote the subtree rooted at $y$, $T_x$ denote the subtree rooted at $x$, excluding those nodes in $T_y$. The remaining subgraph is denoted by $G_{xy}$. 

We consider two cases:  (1) $x$ is not the root of $T$, and (2) $x$ is the root of $T$.

çBy Lemma~\ref{graphDist.lem}, the resistance distance of a node $i$ in $T_x$ (or $T_y$) to the leader set depends only on the resistance distance to $x$ (or $y$). 
Let $R(T_x) = \sum_{i \in T_x} r(i,x)$.
The height of the subtree rooted at $x$  is $h_x=h - d_{xr}$, where $d_{xr}$ is the graph distance between $x$ and $r$.
At each level $i$ in $T_x$ there are $\M^i$ nodes, each at distance $i$ from $x$. Thus,
\begin{equation}
\resizebox{.5 \textwidth}{!}
{
$R(T_x)=\sum_{i=1}^{h_x}\M^i\cdot i=\frac{\M}{(\M-1)^2}\left((\M h_x-h_x-1)\M^{h_x}+1\right)\nonumber.$
}
\end{equation}
A similar expression can be obtained for $R(T_y)$.

We next consider $G_{xy}$.  We can think of this subgraph as a path graph connecting nodes $x$ and $y$, denoted by $P(x,y)$, with each node in the path the root of its own subtree. 
For any node $j$ on the path between $x$ and $y$, $r(j,\{x,y\})$ is given by Lemma~\ref{rpRes.lem}.
For any node $v$ in the subtree $T_j$, its resistance distance to $\{x,y\}$ is
\[
r(v,\{x,y\}) = d_{vj} + r(j,\{x,y\}).
\]
From this, we obtain,
\[
\resizebox{.5 \textwidth}{!}
{
$R(G_{xy}) = \frac{1}{6}({d_{xy}}^2 - 1) + \sum_{\substack{j \in {P(x,y)} \\ j\neq x,y}} \left( R(T_j) + ( |T_j| - 1) r(j,\{x,y\}) \right).$
}
\]
The first term is the total resistance distance for nodes on the path from $x$ to $y$.  For the summation terms,
first we compute the total resistance distance from nodes in the subtree rooted at $j$ to $j$.  Then, for each
node in the subtree, excluding $j$, we add the resistance distance from $j$ to $\{x,y\}$.
An equivalent expression is:
\begin{small}
\begin{align}
R(G_{xy}) =& \sum_{\substack{j \in {P(x,y)} \\ j\neq x,y}} \left( R(T_j) +  |T_j| r(j,\{x,y\}) \right) \nonumber\\
=& \sum_{\substack{j \in {P(x,y)} \\ j\neq x,y}} R(T_j)+\sum_{\substack{j \in {P(x,y)} \\ j\neq x,y}} |T_j| r(j,\{x,y\}). \label{twoSums} 
\end{align}
\end{small}
To simplify the first sum in (\ref{twoSums}), we first consider the subtrees rooted at nodes on the path from $x$ to $r$, denoted by $P(x,r)$ (excluding $r$):
\begin{small}
\begin{align}
\sum_{j\in P(x,r)\atop{j\neq x,r}}R(T_j)=&\frac{d_{xr}-1}{\M-1}+\frac{((\M-1)h-\M-1)\cdot \M^{h}}{(\M-1)^2} \nonumber \\
&+\frac{((\M-1)d_{xr}-(\M-1)h+2)\cdot \M^{h-d_{xr}+1}}{(\M-1)^2}. \label{sum1x.eq}
\end{align}
\end{small}
A similar expression can be obtained for the subtrees rooted at nodes on the path from $r$ to $y$, substituting $d_{xr}$ with $d_{yr}$.

To simplify the second sum in (\ref{twoSums}), we also first consider the subtrees rooted on nodes on the path $P(x,r)$, which is:
\begin{small}
\begin{align}
&\sum_{\substack{j \in {P(x,r)} \\ j\neq x,y}} |T_j| r_j(\{x,y\}) = \M^{h+1}\left(\frac{1}{(\M-1)^2}+\frac{\M+1}{(M-1)^3d_{xy}}\right)\nonumber\\
&~~~+\M^{h}\frac{((\M-1)d_{xr}-\M)((\M-1)(d_{xy}-d_{xr})+\M)}{(\M-1)^3d_{xy}}\nonumber\\
&~~~-\frac{\M^{h+1}}{(\M-1)^3d_{xy}}. \label{sum2x.eq}
\end{align}
\end{small}
As before, a similar expression can be obtained for the subtrees rooted at nodes on the path from $r$ to $y$, substituting $d_{xr}$ with $d_{yr}$.

The above sums (\ref{sum1x.eq}) and (\ref{sum2x.eq}), and their corresponding sums for $P(r,y)$ account for the subtrees rooted at two children of $r$, one containing leader $x$ 
and one containing leader $y$.  For each of the remaining $M-2$ children of $r$, the total resistance distance to $x$ and $y$ from the subtree rooted at child $v$ is
\begin{small}
\begin{equation}
R(T_v)=\frac{((h-1)\M-h)\M^h+\M}{(\M-1)^2}+\frac{\M^h-1}{\M-1} \left(d_{xr}-\frac{d_{xr}^2}{d_{xy}}+1\right).\nonumber
\end{equation}
\end{small}
Combining all of these sums  and including $r(r,\{x,y\})$, we obtain $R(G_{xy})$. Substituting the expressions for $R(T_x)$, $R(T_y)$, $R(G_{xy})$, and the equality $d_{yr} = d_{xy} - d_{xr}$ into (\ref{RofTree.eq}) leads to (\ref{rgkary}).

\noindent \emph{Case 2: $x$ is the root.} In this case, $R(T_y)$ is the same as in Case 1, but $R(T_x)$ is now 
\begin{small}
\begin{align}
R(T_x)=\sum_{i=1}^{h}(\M-1)\M^{i-1}\cdot i=\frac{\M^h\left(\M h-h-1\right)+1}{\M-1}\,.\label{sum3y.eq}
\end{align}
For $R(G_{xy})$, we only need to consider the path from root to $y$ by using (\ref{sum2x.eq}). Combining  (\ref{sum3y.eq}) and (\ref{sum2x.eq}), we obtain
\begin{align}
&\Omega(\{r,y\}) = \frac{d_{xy}}{M-1}+\frac{M}{(M-1)^2}\nonumber\\
&+M^{h+1}\left(\frac{2}{(M-1)^2}+\frac{M+1}{(M-1)^3dy}\right)(M^{-d_{xy}}+1)\nonumber\\
&+M^{h+1}\left(\frac{h}{M-1}-\frac{3}{(M-1)^2}-\frac{2(M+1)}{(M-1)^3d_{xy}}\right)\,,\label{rgkaryXr}
\end{align}
\end{small}
which is equal to (\ref{rgkary}) given $d_{xr}=0$. 

Thus, we conclude that in a perfect $\M$-ary tree, (\ref{rgkary}) holds for any leader set $\{x,y\}$ where their lowest common ancestor is the root.
\end{IEEEproof}

\subsection{Proof of Theorem~\ref{btree.thm}}

\begin{IEEEproof}
We first simplify (\ref{rgkary}) in Lemma~\ref{Mxy.lem} for $M=2$,
\begin{small}
\begin{align}
&\Omega(\{x,y\})  = \textstyle  2 + 2^h(2h -6) + d_{xy} -2^{h+1} \cdot \frac{6}{d_{xy}} \nonumber \\
&~~~~~+ \left(2^{h+1}+1\right)\left( d_{xr} - \frac{d_{xr}^2}{d_{xy}}\right) \nonumber \\
&~~~~~+ 2^{h+1} \left(2^{-d_{xr}} + 2^{-(d_{xy}-d_{xr})}\right) \left( \frac{3}{d_{xy}} + 2 \right). \label{rgbinary}
\end{align}
\end{small}

For a given $d_{xy}$  and $h$, we treat $\Omega$ as a continuous function with argument $d_{xr}$. 
We derive  expressions for  its first and second derivative:
\begin{small}
\begin{align}
\frac{\partial \Omega}{\partial {d_{xr}}}=&\left(2^{h+1}+1\right) \left(1-\frac{2d_{xr}}{d_{xy}}\right)\nonumber\\
& + 2^{h+1} \left(2^{d_{xr}-d_{xy}}-2^{-d_{xr}}\right)\left(\frac{3}{d_{xy}}+2\right)\cdot \ln 2\,, \label{firstOrderD} \\
 \frac{\partial^2 \Omega}{{\partial^2 d_{xr}}} &= \left(2^{h+1}+1\right) \left(-\frac{2}{d_{xy}}\right)\nonumber\\
&\qquad + 2^{h+1}\left(2^{-d_{xr}}+2^{d_{xr}-d_{xy}}\right) \left(\frac{3}{d_{xy}}+2\right) \cdot (\ln 2)^2 \label{secondOrderD}\\
 &\geq \left(2^{h+1}+1\right) \left(-\frac{2}{d_{xy}}\right) + 2^{h+2-\frac{d_{xy}}{2}} \left(\frac{3}{d_{xy}}+2\right)(\ln 2)^2.\label{sndDInE}
\end{align}
\end{small}
From (\ref{firstOrderD}), we observe that $\Omega$ has an extremum at $d_{xr}=d_{xy}/2$.
For $d_{xy} \leq 5$ and $h \geq 4$, (\ref{sndDInE}) is strictly positive, thus $\Omega$ is convex with respect to $d_{xr}$. 
This means that $d_{xr} = d_{xy}/2$ is a minimizer for the given $d_{xy}$.

For  $h\geq 4$, and $d_{xy} \leq 5$, we  examine the potential integer minimizers $d_{xy}=5$, $d_{xr}=2$; $d_{xy}=4$, $d_{xr}=2$; $d_{xy}=3$, $d_{xr}=1$; $d_{xy}=2$, $d_{xr}=1$; and $d_{xy}=1$, $d_{xr}=1$. 
By comparing them in $\Omega(\{x,y\})$ in (\ref{rgbinary}), we  find the minimum is always attained at $d_{xy}=4$, $d_{xy}=2$.

For $d_{xy} \geq 6$ and $h \geq 4$, by checking $\frac{\partial \Omega}{\partial {d_{xr}}}$ and $\frac{\partial^2 \Omega}{{\partial^2 d_{xr}}}$, we  observe that $\Omega$ has two minima. 
Because of the symmetry of the function, these two minima must have the same value, and so we only need to study the solution where $d_{xr}\leq  d_{xy}/2$. 
Since $d_{xy}\geq 6$ and $h\geq 4$, we have:
\begin{small}
\begin{align*}
\frac{\partial \Omega}{\partial {d_{xr}}}\Big|_{d_{xr}=0}=&\left(2^{h+1}+1\right) \\
&+2^{h+1}\Big(\frac{3}{d_{xy}} +2\Big)\ln 2\left(2^{-d_{xy}}-1\right)<0\,,
\end{align*}
\end{small}
and
\begin{small}
\begin{align*}
\frac{\partial \Omega}{\partial {d_{xr}}}\Big|_{d_{xr}=2}=&\left(2^{h+1}+1\right)\left(1-\frac{4}{d_{xy}}\right)\\
&+2^{h+1}\bigg(\frac{3}{d_{xy}} +2\bigg)\ln 2\left(2^{2-d_{xy}}-\frac{1}{4}\right)>0\,.
\end{align*}
\end{small}
Therefore, an integer minimizer of $\Omega$ is attained in the set $d_{xr}\in \{0,1,2\}$.
It is readily verified that for $h \geq 4$ and ${d_{xy} \geq 6}$, 
\[
\Omega|_{d_{xy} \geq 6, d_{xr}=k} \geq \Omega|_{d_{xy}=4,d_{xr}=2}, 
\]
for $k\in\{0,1,2\}$.
This implies that $d_{xy}=4$, $d_{xr}=2$ is  the integer solution that minimizes $\Omega$ for all $h\geq 4$. 
We obtain the expression for $\Rnf$ in (\ref{Rbtree.eq})  by substituting $d_{xy}=4$, $d_{xr}=2$ and $n=2^{h+1}-1$  into (\ref{rgbinary}) and applying $\Omega=2\Rnf$.
\end{IEEEproof}

\subsection{Proof of Theorem \ref{3Theo}}

\begin{IEEEproof}
Based on Lemma \ref{Mxy.lem}, we  derive $\Omega(\{x,y\})$  for a perfect ternary with height $h$, where $x$ and $y$ have the root as their lowest common ancestor,
\begin{small}
\begin{align}
&\Omega(\{x,y\}) = \frac{3^{h+1}+1}{2}\left(d_{xr}-\frac{d^2_{xr}}{d_{xy}}\right)\nonumber\\
&~~~~~~+\frac{3^{h+1}}{2}\left(1+\frac{1}{d_{xy}}\right)(3^{d_{xr}-d_{xy}}+3^{-d_{xr}})\nonumber\\
&~~~~~~+\frac{3^{h+1}}{2}\left(h-\frac{3}{2}-\frac{2}{d_{xy}}\right)+\frac{d_{xy}}{2}+\frac{3}{4}. \label{omega3.eq}
\end{align}
\end{small}
For a given $d_{xy}$, we find its first and second derivative,
\begin{small}
\begin{align}
&\frac{\partial \Omega}{\partial {d_{xr}}}=\frac{3^{h+1}+1}{2} \left(1-\frac{2d_{xr}}{d_{xy}}\right)\nonumber\\
&\quad +\frac{3^{h+1}}{2}\cdot  \left(3^{d_{xr}-d_{xy}}-3^{-d_{xr}}\right)\left(\frac{1}{d_{xy}}+1\right)\cdot \ln 3 \label{firstOrderDforM3} \\
&\frac{\partial^2 \Omega}{\partial^2 {d_{xr}}}=(3^{h+1}+1) \left(-\frac{1}{d_{xy}}\right)\nonumber\\
&\quad +\frac{3^{h+1}}{2} \left(3^{d_{xr}+d_{xy}}+3^{-d_{xr}}\right)\left(\frac{1}{d_{xy}}+1\right)\cdot (\ln 3)^2\nonumber\\
&\geq(3^{h+1}+1) \left(-\frac{1}{d_{xy}}\right)+3^{h+1-\frac{d_{xy}}{2}}\left(\frac{1}{d_{xy}}+1\right)\cdot (\ln 3)^2\,. \label{secondOrderDforM3}
\end{align}
\end{small}
Similar to the proof of Theorem~\ref{btree.thm}, we obtain that $d_{xy}=2$, $d_{xr}=1$ is the optimal integer solution for any $h\geq 4$ and $d_{xy}\leq 2$. By enumerating all $d_{xy}$, $d_{xr}$, given $d_{xy}\leq 4$, we can verify that $d_{xy}=2$, $d_{xr=1}$ is the optimal solution for any $h\geq 4$ and $d_{xy}\leq 4$.

As shown by $\frac{\partial \Omega}{\partial {d_{xr}}}$ and $\frac{\partial^2 \Omega}{{\partial^2 d_{xr}}}$, for a given $d_{xy}$, $\Omega$ has two minima. Because of the symmetry of (\ref{omega3.eq}), we only need to study the minimum that satisfies $d_{xr}\leq d_{xy}$. For any given $d_{xy}\geq 5$, $h\geq 4$, 
\begin{small}
\begin{align*}
&\frac{\partial \Omega}{\partial {d_{xr}}}\Big|_{d_{xr}=0}=\\
&~~~~~\frac{3^{h+1}+1}{2}+\frac{3^{h+1}}{2}\Big(3^{-d_{xy}} -1\Big)\left(\frac{1}{d_{xy}}+1\right)(\ln 3)<0 \\
&\frac{\partial \Omega}{\partial {d_{xr}}}\Big|_{d_{xr}=2}= \\
&\frac{3^{h+1}+1}{2}\left(1-\frac{4}{d_{xy}}\right)+\frac{3^{h+1}}{2}\Big(3^{2-d_{xy}} -\frac{1}{9}\Big)(\frac{1}{d_{xy}}+1)(\ln 3)>0.
\end{align*}
\end{small}
Thus, the optimal real-valued $d_{xr}$ lies in $(0,2)$. By evaluating (\ref{omega3.eq}) for $k\in\{0,1,2\}$, it can be verified that,
\[
\Omega|_{d_{xy} \geq 5, d_{xr}=k} \geq \Omega|_{d_{xy}=2,d_{xr}=1}.
\]
Thus, we have shown that $d_{xy}=2$, $d_{xr}=1$ is the global integer minimizer for all $h\geq 4$ in perfect ternary trees. We obtain (\ref{Rttree.eq}) by substituting $d_{xr}=1$, $d_{xy}=2$ into (\ref{rgbinary}) and applying $n=(3^{h+1}-1)/2$ and $\Omega=2\Rnf$.
\end{IEEEproof}

\subsection{Proof of Theorem~\ref{KTheo}}


\begin{IEEEproof}
We start by calculating $\frac{\partial \Omega}{\partial d_{xr}}$ and $\frac{\partial^2 \Omega}{\partial^2 d_{xr}}$.,
\begin{small}
\begin{align}
\frac{\partial \Omega}{\partial d_{xr}} &= \frac{\M^{h+1}+1}{\M-1}\left(1-\frac{2d_{xr}}{d_{xy}}\right) +(\ln M)M^{h+1}\Big(\frac{2}{(M-1)^2} \nonumber \\
&~~~~ +\frac{M+1}{(M-1)^3dy}\Big)(M^{d_{xr}-d_{xy}}-M^{-d_{xr}}) \label{fdM.eq} \\
\frac{\partial^2 \Omega}{\partial^2 d_{xr}} &= \frac{\M^{h+1}+1}{\M-1}\left(1-\frac{2}{d_{xy}}\right)+(\ln M)^2 M^{h+1}\Big(\frac{2}{(M-1)^2} \nonumber \\
&~~~~ +\frac{M+1}{(M-1)^3dy}\Big)(M^{d_{xr}-d_{xy}}+M^{-d_{xr}}). \label{sdM.eq}
\end{align}
\end{small}
From (\ref{fdM.eq}) and (\ref{sdM.eq}), we observe that $\Omega$ has a minimum that satisfies $d_{xr}\leq \frac{d_{xy}}{2}$. 
Since $\Omega$ is symmetric about $d_{xr}=\frac{d_{xy}}{2}$, we only consider potential integer minimizers with $d_{xr} \leq \frac{d_{xy}}{2}$. 
Further,
\begin{small}
\begin{align}
&\frac{\partial \Omega}{\partial d_{xr}}\Big|_{d_{xr}=1}=\frac{M^h}{(M-1)^3}\bigg(\left(M+\frac{1}{M^h}\right))(M-1)^2(1-\frac{2}{d_{xy}}) \nonumber \\
&\qquad +\left(2(M-1)+\frac{M+1}{d_{xy}})(M^{2-d_{xy}}-1\right)\ln(M)\bigg), \label{fdMxr.eq}
\end{align}
\end{small}
and observe that $\frac{\partial \Omega}{\partial d_{xr}}|_{d_{xy}=2,d_{xr}=1}=0$. 
For $d_{xy} \geq 3$, we can lower bound (\ref{fdMxr.eq}) by  
\begin{small}
\begin{align}
&\frac{\partial \Omega}{\partial d_{xr}}\Big|_{d_{xr}=1} >\frac{M^h}{(M-1)^3}\bigg(M(M-1)^2(1-\frac{2}{d_{xy}}) \nonumber \\
&\qquad -(2(M-1)+\frac{M+1}{d_{xy}})\ln(M)\bigg). \label{fbound.eq}
\end{align}
\end{small}
For $d_{xy} \geq 3$, the bound (\ref{fbound.eq}) is positive for $M=4$ and $h\geq 4$, and it is increasing in $M$, $d_{xy}$ and $h$.
Thus, for all $d_{xy} \geq 2$, the integer minimizer of $d_{xr}$ will be either 0 or 1.
Further, for $d_{xy} = 1$, the only potential solution that satisfies $d_{xr}\leq d_{xy}/2$ is $d_{xr} =0$.

%

We propose that the optimal integer solution is ${d_{xy} =1}$, ${d_{xr} = 0}$, and we validate its optimality by comparing with 
$d_{xy} \geq 2$ and $d_{xr} \in \{0,1\}$.
For $d_{xy} \geq 2$,
\begin{small}
\begin{align*}
&\Omega|_{d_{xr}=0}-\Omega|_{d_{xy}=1,d_{xr}=0} =\\
&~~~~~~~~\frac{1}{(M-1)^3d_{xy}}\Big((M-1)^2d_{xy}(d_{xy}-1)\\
&~~~~~~~~\quad +M^{h+1-d_{xy}}((2d_{xy}+1)M-2d_{xy}+1)\\
&~~~~~~~~\quad +M^{h}(d_{xy}(M-1)^2-M(M+1))\Big)\,,
\end{align*}
\end{small}
which is positive for $M\geq 4$. In addtion,
\begin{small}
\begin{align*}
&\Omega|_{d_{xr}=1}-\Omega|_{d_{xy}=1,d_{xr}=0}=\\
&~~~~~~~\frac{1}{(M-1)^3d_{xy}}\Big((M-1)^2(d^2_{xy}-1)\\
&~~~~~~~\quad +M^{h+2-d_{xy}}((2d_{xy}+1)M-2d_{xy}+1)\\
&~~~~~~~\quad +M^{h}(d_{xy}(M-1)^3-M(M^2+2))\Big)\,,
\end{align*}
\end{small}
 is also positive for $d_{xy}\geq 2$ and $M\geq 4$. 
Therefore, the optimal leader set $\{x,y\}$ is such that $d_{xy}=1$, $d_{xr}=0$ when $M\geq 4$ and $h\geq 4$. Then, (\ref{Rktree.eq}) is obtained by substituting $d_{xr}=0$, $d_{xy}=1$, $n=(\M^{h+1}-1)/(M-1)$ into (\ref{rgkary}) 
and using the fact that $\Omega=2\Rnf$.
%
\end{IEEEproof}

\bibliographystyle{IEEEtran}
\bibliography{leader}

%

%
%
%





\end{document}